\documentstyle{amsppt}

\def\id{\operatorname{id}}
\def\alg{\operatorname{alg}}
\def\II{\operatorname{II}}
\def\Span{\operatorname{Span}}

\topmatter

\title  Nonfree actions of countable groups and their characters
\endtitle

\author A.~M.~Vershik \endauthor

\thanks
Partially supported by the RFBR grants
08-01-00379-a and
09-01-12175-ofi-m.\endthanks



\address St.~Petersburg Department of Steklov Institute of Mathematics,
St.~Petersburg, Russia.
\endaddress
\email vershik\@pdmi.ras.ru
\endemail

\abstract
We introduce a number of definitions of nonfree actions of groups. The
most important of them is that of a totally nonfree action; it is
naturally related to the theory of characters of groups and their
factor representations. This short note is a brief exposition of a
part of a more detailed paper on this subject, which is now in
preparation.
\endabstract

\endtopmatter
\document

\hfill  {\bf To Rita's memory }

\head 1.  Definitions  \endhead

Let $G$ be a countable group acting on a Lebesgue space $(X, {\goth
A}, \mu)$ with
continuous measure. Given an element $g \in G$, denote by
$G_X$ the (measurable) set of fixed points of $g$:
$G_x=\{x\,:\,gx=x\}$. According to a well-known definition, the action is called free if
 $\mu X_g=0$ for  $g \ne \id$.
Denote by ${\goth A}_G$ the $\sigma$-subalgebra of
$\goth A$ generated by all such sets of fixed points:
${\goth A}_G=\sigma\text{-}\alg \{X_g; g\in F\}$.
For a free action, this $\sigma$-algebra is trivial.

\proclaim{Definition 1}
The action of $G$ is called totally nonfree if the $\sigma$-algebra
${\goth A}_G$ coincides with the whole original $\sigma$-algebra:
${\goth A}_G=\goth A$.
\endproclaim

We can give a similar definition for arbitrary (uncountable) groups,
but in this case we should assume that the action is individual
(i.\,e., there is a set of full measure for all points of which the
action of all group elements is well defined).

Consider the {\it lattice} (with respect to the natural order)
$L(G)$ {\it of all subgroups of $G$}. This lattice has been
intensively studied from the purely group-theoretic point of view
since the 1930s,
and even earlier in the 19th century (see the recent monograph
\cite{2}).   In a
natural way we can introduce a weak topology on the lattice $L(G)$ (a neighborhood of a given subgroup $H$ is the family of all
subgroups for which the set of words of length at most $n$
coincides with the same family for $H$) and the
corresponding Borel structure. Obviously, in this topology the lattice $L(G)$ is
compact and is a totally disconnected
uncountable separable Cantor set. The group $G$ acts on $L(G)$
continuously by conjugation; this action is sometimes called {\it
adjoint}. We will study Borel conjugation-invariant probability
measures on lattices of subgroups. Apparently, problems concerning
these measures and the dynamics of the action of $G$ on $L(G)$ have not
yet been studied.

With almost every point $x \in X$ we associate its stabilizer, i.e.,
the subgroup $G_x=\{g\,:\, gx=x\}$, regarded as an element of the
lattice, and denote this map by
$\Psi_G:X\rightarrow L(G)$:
$$
\Psi_G: x\mapsto G_x.
$$

Obviously, the map $\Psi_G$ is Borel-measurable.

Denote the $\sigma$-algebra generated by this map (i.e., the inverse
image of the  $\sigma$-algebra of Borel sets on
$L(G)$) by ${\goth A}^G$.

\proclaim{Lemma 2}
$$
{\goth A}_G \subseteq {\goth A}^G,
$$
the equality takes place for countable groups.
\endproclaim

Indeed, if the  $\sigma$-algebra ${\goth A}_G$ is trivial, then
the inclusion is proved. Assume that
${\frak A}_G$ is nontrivial, and let $\xi$ be the measurable partition
(equivalence relation) corresponding to the $\sigma$-algebra ${\goth
A}_G$. Obviously, by the definition of a measurable partition, two
points belonging to the same block of $\xi$ have the same stabilizer
and hence lie in the inverse image under
$\Psi_G$ of the same point.

\proclaim{Corollary 3}
If the action of $G$ is totally nonfree, then the map
$\Psi_G$ is a monomorphism $\bmod 0$; in other words, the
$\sigma$-algebra ${\goth A}_G$ coincides with
${\goth A}^G$.
\endproclaim

\proclaim{Definition 4}
Let us say that the action of $G$ is extremely nonfree if the map
$\Psi_G$ is monomorphic; in other words, if the $\sigma$-algebra
${\goth A}^G$ coincides with the whole $\sigma$-algebra
$\goth A$. In terms of the action, this means that almost all points
have different stabilizers.
\endproclaim

It is easy to give an example of uncountable group
distinguishing between the total and extremal nonfreeness of
actions. For instance (though the group here is continuous),
the action of the group $SO(3)$ on the projective plane $P_2R$
is extremely nonfree, but not totally nonfree.

Consider the image  $\nu_{\mu}$ of the measure $\mu$ on
$L(G)$ under the map
$\Psi_G$. Obviously, $\nu_{\mu}$ is invariant under conjugation.
Observe the following simple but important fact.

\proclaim{Theorem 5}
In the class of extremely nonfree actions of the group $G$, the measure
$\nu_{\mu}$ defined above is a complete metric invariant of the
action. In other words, two extremely nonfree actions of the group
$G$  in spaces
$(X,{\goth  A},\mu)$ and $(X',{\goth A'},\mu')$ are metrically
isomorphic if and only if the measures
$\nu_{\mu}$ and $\nu_{\mu'}$ on the lattice $L(G)$ coincide.
\endproclaim

This follows immediately from the fact that the map
$\Psi_G$, which transfers isomorphically an extremely nonfree action from
an arbitrary Lebesgue space to the lattice
$L(G)$, is monomorphic.

\head 2. Extremely nonfree actions \endhead

We will consider conjugation-invariant measures on the lattice of
subgroups. Examples of such measures are the delta-measure at the
identity or the atomic measure on the set of conjugate subgroups with
normalizer of finite index. Moreover, for every action of the group
$G$ with an invariant measure $\mu$ on an arbitrary space $(X, {\goth
A},\mu)$, the map  $\Psi_G$ sends $\mu$ to a conjugation-invariant
measure on the lattice $L(G)$. We are interested only in continuous
(and nontransitive in the case of continuous groups) ergodic
probability measures.

\proclaim{Problem 1} For what countable groups do there exist
conjugation-invariant continuous ergodic
probability measures? Describe all such measures.
\endproclaim

Examples of such groups are the infinite symmetric and similar groups;
the list of all such measures for them is known (see below). It is
recently shown in \cite{4}, as an answer to the above question, that for
noncommutative free groups and similar groups, such measures do exist;
however, we still
do not have their complete description.

Let $\nu$ be such a measure on the lattice $L(G)$. Is it true that the
action of the group $G$ on the space
$(L(G), \nu)$ is extremely nonfree? In general, this is not the case;
moreover, simple examples show that even if the action of $G$ on the
space $(X,\mu)$ is extremely nonfree, in general, the
$\Psi_G$-image $\nu_{\mu}$ of the measure $\mu$ is such
that the action of $G$ on the space $(L(G),\nu_{\mu})$ is not extremely nonfree. The reason is as
follows: since the stabilizer of a subgroup $H$, regarded as an
element of the lattice $L(G)$, is its normalizer
$N(H)=\{g\in G: gHg^{-1}=H \}$, it is quite possible that the
normalizers coincide for different subgroups. If the action of the
group $G$ on the lattice
$L(G)$ with some measure
 $\nu$ by conjugation is extremely nonfree, we will say that the measure
$\nu$ itself is {\it extremely nonfree}.

\proclaim{Definition 6}
A subgroup $H$ of an arbitrary group $G$ is called
abnormal\footnote{The author does not know whether there is a
generally accepted term for such subgroups.} if it coincides with its
normalizer.
\endproclaim

\proclaim{Proposition 7}
A conjugation-invariant measure on
$L(G)$ is extremely nonfree if and only if the measure of the set of
abnormal subgroups is equal to one.
\endproclaim

\demo{Proof}
Assume that a measure $\nu$ is extremely nonfree, but there is a set
of positive $\nu$-measure that consists of abnormal subgroups, i.e.,
subgroups $H$ such that for each of them we can find an element
$h_H\in N(H)\setminus H$. Since $G$ is countable, there exists an
element $h$ such that the last condition holds for all subgroups from
a set of positive measure. Then for a subgroup $H$ from this set we have
$hHh^{-1}\ne H$, but at the same time the normalizers of both
subgroups, which are stabilizers for the adjoint action, coincide,
contradicting the extremal nonfreeness.
The converse is obvious.
\hfill $\square$
\enddemo

Consider the normalization operation
$\Cal N$ on the lattice of subgroups, which associates with every subgroup its normalizer:
${\Cal N}(H)=N(H)$. It can be extended to an operation on measures:
${\Cal N}_*(\nu)(E)=\nu({\Cal N}(E))$. Obviously, if a measure
$\nu$ is invariant, then its image is also invariant. If
the measure ${\Cal N}_*(\nu)$ is extremely nonfree, then we say that
the measure $\nu$ is reducely extremely nonfree; this case occurs
often and hence should be distinguished. As follows from the above
proposition, in terms of the group this means that the measure $\nu$ is
concentrated on subgroups $H$ satisfying the condition
$$
N(N(H)=N(H)),\quad \text{or} \quad  {\Cal N}^2(H)={\Cal N}(H).
$$

However, this is not always the case, and the normalization process
does not in general stabilize after the first step, or any finite number of
steps.\footnote{The author does not know corresponding examples, but,
apparently, stabilization may not be achieved after infinite and even
countable number of steps.}

Anyway, we are mainly interested in extremely nonfree measures on
lattices.

\proclaim{Problem 2} For what countable groups do there exist continuous
ergodic probability extremely nonfree measures invariant under
conjugation?
\endproclaim

Apparently, the class of such groups is significantly more narrow than
the class of groups solving the first problem.

\head 3. Totally nonfree actions\endhead

Totally nonfree actions are defined above as actions for which the
collection
of all sets of fixed points of different group elements generates the
whole $\sigma$-algebra; as we have seen, these actions are extremely
nonfree, so that we can study them on lattices of subgroups.
Consider an extremely
nonfree measure $\nu$ on the lattice $L(G)$ of subgroups of a group
$G$. This measure is concentrated on the set of abnormal subgroups;
let us formulate the condition of total nonfreeness. Following the
tradition, outlined above, to extend the term describing an action to the
corresponding measures, {\it we say that a measure
$\nu$ on the lattice $L(G)$ is totally nonfree if the adjoint action
of $G$ on the measure space
$(L(G),\nu)$ is totally nonfree}. Given an element
$g\in G$, we denote by
$L_g\subset L(G)$ the set of subgroups that contain $g$.

\proclaim{Proposition 8}
An extremely nonfree measure on the lattice of subgroups $L(G)$
is totally nonfree if and only if the collection of sets
$L_g$, $g\in G$, that are of positive
$\nu$-measure generates the whole $\sigma$-algebra in the space
$(L(G),\nu)$; in other words, $\nu$-almost every pair of different
subgroups can be distinguished from each other by a set
$L_g$ of positive $\nu$-measure.
\endproclaim

Indeed, the set of fixed points for an element $g\in G$ is the set of
subgroups whose normalizers contain $g$; but for an extremely
nonfree measure, almost all subgroups are abnormal, i.e., coincide
with their normalizers, hence this is the set of subgroups containing
$g$.

The most interesting problem involves totally nonfree actions.

\proclaim{Problem 3} For what countable groups do there exist continuous
ergodic totally nonfree probability measures
invariant under conjugation? Describe all such measures. In a
more general form, for what groups do there exist totally nonfree
actions? Describe all such actions up to isomorphism.
\endproclaim

As it will become clear from the next section, it is this class of
groups that is important for representation theory. Examples are the
infinite symmetric group
${\goth S}_{\Bbb N}$, the group $GL(\infty ,q)$ of infinite matrices
over a finite field, etc.

\head 4. Relation to representations and the theory of characters \endhead

Every action of a group $G$ on a space
$(X,\mu)$  with invariant measure generates a unitary representation of $G$ in the space
$L^2_{\mu}(X)$ according to the formula
$[U_g(f)](x)=f(gx)$ (Koopman representation). The question about the
properties of this representation and its irreducible
decomposition is difficult even for
$\Bbb Z$ (this is the subject of the spectral theory of dynamical
systems). Since the measure is finite, we have the invariant
one-dimensional subspace of constants. But the following problem
remains open.

\proclaim{Problem 4} In what cases is the Koopman representation irreducible in the
orthogonal complement to the subspace of constants?
\endproclaim

Obviously, the ergodicity of the action is a necessary condition,
because the decomposition into ergodic components gives rise to the
decomposition of the Hilbert space
$L^2_{\mu}(X)$ into a direct integral (or direct sum) of invariant
subspaces.

The problem can be formulated  more specifically. Consider an ergodic
action of a group $G$ on a space
$(X,\mu)$ with invariant measure and the
$W^*$-algebra generated by the operators of multiplication
by bounded measurable functions and the operators of the action of $G$. This
algebra may be regarded as the image of the skew product
$ l^1(G) \rightthreetimes L^{\infty}(X)$
of the group algebra of $G$ and the algebra
$L^{\infty}(X)$ of multiplication operators. It turns out that it coincides with the algebra of
all bounded operators; indeed, the algebra of multiplication operators
is a maximal commutative self-conjugate subalgebra in the algebra of
all bounded operators, and, besides, by ergodicity, it has no
nonconstant multiplicators
commuting with the group action, so that the commutant of this
algebra is scalar, and our assertion follows from the von Neumann
bicommutant theorem. Thus Problem 4 can be equivalently formulated  as
follows.

\proclaim{Problem 4$'$}
In what cases does the $W^*$-algebra spanned by the operators of the group
action contain all multiplicators with zero integral?
\endproclaim

Most likely, Problem 4 is very difficult. Examples of a positive
answer the author is aware of are quite rare.

There exists another canonical representation associated with a
measure-pre\-serv\-ing action of a group, which goes back to von Neumann
(see, e.g., \cite{3}). It can be called the {\it trajectory}, or {\it
groupoid} representation; see, e.g., \cite{5}. Consider the graph $\Pi$
of the group action, i.e., the set of pairs
$\{(x,y): y=gx, \, g\in G\}$, regarded as a measurable partition in
$X\times X$. Endow it with the $\sigma$-finite measure
$M$ that induces the
measure  $\mu$ on both factors  $X\times *$ and $*\times X$ and whose
conditional measure in the layer over each point
$(*,y)$ or $(x,*)$ is the
 counting (i.e., uniform infinite) measure for all $x$ and $y$. Then on
$\Pi$ we have two commuting actions of the group $G$ with invariant
measure $M$ --- the left one,
$(x,y)\mapsto (gx,y)$,  and the right one, $(x,y)\mapsto(x,gy)$.
Correspondingly, in the space
$L^2_M(\Pi)$ we have two unitary representations of the group
$G$ and two $*$-representations of the skew product mentioned above.
It is well known that if the action of $G$ on the space
$(X,\mu)$ is ergodic, then both representations of the skew product
are factor representations of type
$\II_1$, the left and the right factors being mutual commutants.
The characteristic function of the diagonal
${\Delta}=\{(x,x),x \in X\}$,
i.e., the element $1_{\Delta} \in L^2_M(\Pi)$, is a bicyclic vector
for the factors. This construction of a representation of the skew
product is called the von Neumann, or groupoid, or trajectory (since
layers are trajectories of the group action) construction. For a free
action, it was suggested and studied in the pioneering papers by von
Neumann; nonfree actions were considered, e.g., by Krieger.
But the case of an extremely nonfree action first appeared in
\cite{7}, where it was considered as a particular example related to
factor representations of the infinite symmetric group and its relatives.

As in the previous example, we can restrict this representation to the
group itself and ask a similar
question about these representations.

\proclaim{Problem 5} In what cases does the restriction of the groupoid
representation to the group generate the whole factor? Or,
equivalently, when does the $W^*$-algebra spanned by the group operators
contain all multiplicators and thus coincide with the whole factor?
\endproclaim

Recall that a complex-valued function
$\phi$ on a group is called a character if it is nonnegative definite,
central, and normalized:
$$
\{\phi (g_ig_j^{-1})\}_{i,j=1}^n \geq 0 ,\quad
\phi(ghg^{-1})=\phi(h), \quad \phi({\id})=1.
$$
The characters form a convex compact set in the weak topology, and
extreme points of this set are called indecomposable characters.

The importance of Problem 5 can be seen from the following
correspondence between the indecomposability of a character and the
factorness.

\proclaim{Theorem 9}
Consider an ergodic action of a group $G$ on a space
$(X,\mu)$ with invariant measure; the function
$$
\phi(g)=\mu(X_g)=\mu\{x: gx=x, x \in X \} \qquad (*)
$$
is a character of $G$. If the restriction of the groupoid
representation to the group generates the whole factor, then this
character is indecomposable.
\endproclaim

\demo{Proof}  First of all, we have the formula
$$
<U_g 1_{\Delta},1_{\Delta}>=\mu(X_g),
$$
which implies the first assertion. As we know from the theory of von
Neumann algebras, a trace on a $W^*$-algebra is indecomposable
(i.e., cannot be written as a convex combination of
other traces) if and only if the corresponding representation is a
factor. In our case, the trace corresponds exactly to a group character.
\hfill  $\square$
\enddemo

Thus we now can search for characters using group
actions. We will see that sometimes this method, with some
supplements, allows one to describe all characters of a group.

In contrast to Problem 4, Problem 5 turns out to have a very lucid
solution, which is precisely the main result of this paper.

\proclaim{Theorem 10}
The character defined by formula
$(*)$ is indecomposable, and thus generates a factor representation of
the group $G$ of type $\II_1$, if and only if the action is totally
nonfree.
\endproclaim

\demo{Proof}
We present a sketch of the proof. For definiteness, consider the left
factor representation of the group $G$ in the space
$L^2_M(\Pi)$, where the set
$\Pi \subset X\times X$ and the measure $M$ were constructed above. We
will prove that the cyclic hull of the characteristic function of the
diagonal $1_{\Delta}$ with respect to the left action of the group
coincides with the whole space
$L^2_M(\Pi)$ if and only if the action of $G$ on the space
$(X, \mu)$ is totally nonfree.

First consider the ``diagonal subspace'' $H_{\Delta}$ of functions from
$L^2_M(\Pi)$ supported by the diagonal
${\Delta}=\{(x,x),x \in X\}$, and let $R_{\Delta}\equiv R$ be the
orthogonal projection to
$H_{\Delta}$. Let us show that if, and only if, the group action is totally
nonfree, then
the linear hull of the projections of the images of the vector
$1_{\Delta}$ under the action of $G$, i.e.,
$\Span \{ R[g (1_{\Delta})]; g \in G\}$, is everywhere dense in the
subspace $H_{\Delta}$. Note that the projection  $ R [g
(1_{\Delta})]$, regarded as a function on the diagonal, is precisely
the indicator function of the set $X_g$ of fixed points
of $g$.  The linear hull of all these indicator functions
coincides with the linear hull of the indicators of subsets from the
subalgebra generated by
the subsets of fixed points of all elements of $G$. Indeed, this
linear hull contains the indicators of the complements of subsets of
fixed points,  since $1_{\bar X_g}=1-1_{X_g}$; it
also contains the indicators of intersections of such subsets, since the
intersection  $X_g \cap X_h$ coincides with
$X_{gh}$ provided that $gx \ne h^{-1}x$ almost everywhere, which we
can assume without loss of generality. Thus the indicator of the
union of two subsets of fixed points also lies in this linear hull,
and hence it contains the indicators of all sets from the
$\sigma$-subalgebra generated by the sets of fixed
points. If, and only if, the action is totally nonfree,
then, by definition,
this $\sigma$-subalgebra is dense in the $\sigma$-algebra of all
measurable sets, and thus the linear hull
$\Span \{R [g (1_{\Delta})]; g \in G\}$ is everywhere dense in the
subspace $H_{\Delta}$. The image of $1_{\Delta}$ can
be written as a sum of two terms:
$g [(1_{\Delta})]=R  [g (1_{\Delta})]\dotplus (I -
R_{\Delta}) [g (1_{\Delta})]$. The support of the second term lies
outside the diagonal, and, applying the action of group elements, we
can ``drive it to infinity,'' i.e., we can shift the support
of $(x,gx)$ by changing $g$. Meanwhile, the first term will still
belong to the subspace $H_{\Delta}$, and the second term will weakly
tend to zero. As a result, in the limit we obtain an element from the linear hull mentioned above supported
by the diagonal. But bounded
functions supported by the diagonal in the left (and also right)
representation act as multiplicators. Thus we have proved that if, and
only if, the action is totally nonfree, then the weak
closure of the algebra of operators of the group action contains all
multiplicators, and this means precisely that the representation of
the group generates the whole factor formed by the skew product of the
group algebra and the Abelian group of multiplicators.
\hfill $\square$
\enddemo

Thus the notion of a totally nonfree action allows one to solve
Problem~5, but does not provide a solution to Problem~4: there exist
examples of totally nonfree actions of two groups such that for one of
them the answer is
positive, and for the other one, the answer is negative. The reason is
as follows: in order to obtain a positive solution to Problem~5, we need to prove
only that the single vector $1_{\Delta}$
is cyclic, and in order to obtain a positive solution to Problem~4, we need to prove that every vector that is
not a constant and is not orthogonal to the subspace of constants, is cyclic,
which is a more difficult task.

Along with Theorem 9, Theorem 10 provides a method for finding
indecomposable characters of groups: one should search for totally
nonfree actions of the group, and the measure of the set of fixed
points of a given element for such an action is the value of the
corresponding character at this
element. This construction has another reserve component: if the group
has
cocycles with values in the circle $\Bbb R/\Bbb Z$ that are not
cohomologous to 1, then, using them in the framework described
above as multiplicators, one can obtain representations nonequivalent to
the previous ones and new characters.
The following problem is important.

\proclaim{Problem 6}  For what groups does this construction give all
indecomposable characters of the group?
\endproclaim

In turns out that this is the case for the infinite symmetric group,
as follows from a new interpretation of the papers by S.~V.~Kerov and
the author published in the 1980s
(\cite{7} and others).

By the way, observe (this concerns the construction of factor
representations) that considering the case of equal frequencies
in Thoma's theorem separately, as in
\cite{7} and subsequent papers, turned out to be superfluous. The
reason is that the old construction used the groupoid of sequences
with a Bernoulli measure, and the present construction is based on the
groupoid of subgroups; hence, whatever the probabilities are, applying
Theorem~9, we obtain the whole factor and not a proper subfactor as in
the case of equal frequencies in the paper \cite{7}.

It seems that the situation is similar for the groups
$GLB(q)$, $GL(\infty,
F_q)$, the group of rational rearrangements of the interval, etc. Note
also that our considerations apply to the infinite-dimensional groups
$U(\infty)$, $O(\infty)$, to the so-called Dye's trajectory group
\cite{8}, etc.; for continuous groups, one should slightly modify some
of the definitions and arguments given above. The transition to the
groupoid of conjugate subgroups makes the construction of
representations completely invariant. Besides, this groupoid is of
interest in itself.

In the paper \cite{1}, which is under preparation, we describe all
extremely nonfree and totally nonfree actions of the infinite
symmetric group, i.e., essentially, the corresponding measures on the
lattice of subgroups of this group, and show how one can obtain from
this description the complete list of characters of this group.
This is another, ``dynamical'' (i.e., based on considering group actions)
proof of Thoma's theorem about the characters of the infinite
symmetric group. A similar proof is
given for closely related groups. It turns out that the list of totally nonfree
measures on the lattice of subgroups is exhausted by the Bernoulli
measures on Young subgroups with infinite blocks; this result is yet another
analog of the famous de Finetti theorem.

In conclusion, we would like to mention that the classes of actions of countable
groups introduced in this paper, and groups for which such actions and
the corresponding measures on the lattices of subgroups do exist, are
of interest beyond the theory of representations and characters.
Apparently, such groups are distinguished also by other properties.

\bigskip

Translated by N.~V.~Tsilevich.

\Refs

\item{1.}
{
A.~Vershik,
{\it Non-free actions of countable groups and their representations}, in preparation.
}

\item{2.}
{
R.~Schmidt,
{\it Subgroup Lattices of Groups}. Walter de Gruyter, Berlin--New York,
1994.
}

\item{3.}
{
J.~Dixmier,
{\it Les $C^*$-alg\`ebres et leur representations}. Gauthier-Villars, Paris, 1964.
}

\item{4.}
{
D.~D'ngeli, A.~Donno, M.~Matter,  and T.~Nagnibeda,
{\it Schreier graphs of the Basilica group}. --- J. Modern Dynamics {\bf 4}, no.~1 (2010), 167--205.
}

\item{5.}
{
M.~Takesaki,
{\it Theory of Operator Algebras}, Vol. 3. Springer-Verlag, Berlin, 2003.
}

\item{6.}
{
J.~Renault,
{\it Groupoid Approach to C*-Algebras}.   Lect. Notes in Math. {\bf793},
Springer-Verlag, Berlin--New York, 1980.
}

\item{7.}
A.~M.~Vershik and S.~V.~Kerov.
{\it Characters and factor representations of the infinite symmetric
group.}
---  Dokl. Akad. Nauk SSSR  {\bf 257}, No.~5 (1981), 1037--1040.

\item{8.}
{
H.~Dye,
{\it On groups of measure preserving transformations}. I, II. ---
Amer. J. Math. {\bf 81}, No. 1 (1959), 119--159; {\bf 85}, No. 4 (1963), 551--576.
}

\endRefs

\enddocument